\documentclass{article}

\usepackage[
  lang = american, 
 paper = hardcover, openaccess,
]{ems-ecm}

\numberwithin{equation}{section}
\newcommand{\snr}[1]{\lvert #1\rvert}
\newcommand{\nr}[1]{\lVert #1 \rVert}
\newtheorem{definition}{Definition}
\newtheorem{remark}{Remark}
\newtheorem{cor}{Corollary}
\newtheorem{theorem}{Theorem}
\newcommand{\mf}[1]{\mathfrak{#1}}
\def\aaa{\mf{a}}
\newcommand{\medint}{-\kern -,375cm\int}
\newcommand{\medintinrigo}{-\kern -,315cm\int}
\def\dx{\,{\rm d}x}
\def\loc{\operatorname{loc}}
\usepackage{tikz}
\usepackage{graphicx}
\usetikzlibrary{patterns}

\newlength{\defbaselineskip}
\newcommand{\mint}{\mathop{\int\hskip -1,05em -\, \!\!\!}\nolimits}
\def\eqn#1$$#2$${\begin{equation}\label#1#2\end{equation}}
\def\dist{\,{\rm dist}}

\numberwithin{equation}{section}

\begin{document}


\title{Sketches of Nonuniformly Elliptic Schauder Theory}
\titlemark{Sketches of Nonuniformly Elliptic Schauder Theory}



\emsauthor{1}{
	\givenname{Cristiana}
	\surname{De Filippis}
	\mrid{1226897} 
	\orcid{0000-0002-3432-1908}}{C.~De~Filippis} 
	
\Emsaffil{1}{
	\department{Dipartimento SMFI}
	\organisation{Università di Parma}
	\address{Viale delle Scienze 53/a, Campus}
	\zip{43124}
	\city{Parma}
	\country{Italy}
	\affemail{cristiana.defilippis@unipr.it}}
%



\begin{abstract}
Schauder theory is a basic tool in the study of elliptic and parabolic PDEs, asserting that solutions inherit the regularity of the coefficients. It plays a central role in establishing higher regularity for solutions to a broad class of elliptic problems exhibiting ellipticity, including those involving free boundaries. In the linear setting, Schauder theory dates back to the 1920-30s and is now considered classical. Nonlinear extensions were developed in the 1980s. All these classical results are restricted to uniformly elliptic operators and heavily rely on perturbative techniques—freezing the coefficients and comparing the solution to that of a constant-coefficient problem. However, such methods fail in the nonuniformly elliptic setting, where homogeneous a priori estimates break down and standard iteration arguments no longer apply. Here we give a brief survey on recent progresses including the solution to the longstanding
problem of proving the validity of Schauder estimates in the nonlinear, nonuniformly elliptic setting.

\end{abstract}

\maketitle

\noindent\textbf{Notation.}~Through this note, $\Omega\subseteq \mathbb{R}^n$, $n\ge 2$, denotes an open, bounded domain with Lipschitz boundary; ${\rm B}\equiv B_{\varrho}\equiv B_{\varrho}(x_{0})\Subset \Omega$ will denote a ball compactly contained into $\Omega$, centered at $x_{0}$ with radius $\varrho>0$. We shall also denote $\mathcal B_\varrho =B_\varrho(0_{\mathbb R^n})$. 
Sometimes we shall use symbols "$\gtrsim$", "$\lesssim$" with subscripts, to indicate that a certain inequality holds up to constants whose dependencies are marked in the subfix. With $\Omega_0 \subset \mathbb{R}^{n}$ being a measurable subset with bounded positive Lebesgue measure $0<| \Omega_0|<\infty$, $\mathbb{N}\ni k\ge 1$, $\kappa\in \mathbb{R}$ being numbers and $w \colon \Omega_0 \to \mathbb{R}^{k}$, being an integrable map, we denote $$(w)_{\Omega_0}\equiv \mint_{\Omega_0}w(x)\dx:=\frac1{\snr{\Omega_0}}\int_{\Omega_0}w(x)\dx$$ its integral average and $(w-\kappa)_{+}:=\max\{w-\kappa,0\}$ its upper truncation. Finally, whenever introducing various objects instrumental to our presentation, unless clearly stated, we will implicitly assume that all the quantities involved are regular enough to ensure the exposition is well-defined.

\section{Measuring ellipticity: equations and functionals}
The Laplace equation $\Delta u=0$ is the most prominent among elliptic partial differential equations of the form
\begin{equation}\label{1.0}
-\textnormal{div} \ A(x,Du)=0\quad \mbox{in} \ \ \Omega,
\end{equation}
defined by means of a Carath\'eodory regular\footnote{This means that $x\mapsto A(x,z)$ is measurable for every choice of $z\in \mathbb R^n$ and $z\mapsto A(x,z)$ is continuous  for a.e. $x\in \Omega$. It serves to ensure that the composition $x\mapsto A(x,D(x))$ is measurable whenever $D$ is a measurable vector field.} vector field $A\colon \Omega \times \mathbb R^n \to \mathbb R^n$. Ellipticity means that $\partial_{z} A(x,z)$, whenever exists, is a (symmetric), positive definite $\mathbb R^{n\times n}$-matrix\footnote{For simplicity we assume that $\partial_{z} A$ is symmetric, i.e. $\partial_{z_j}A_i=\partial_{z_i}A_j$ for every $1\leq i,j\leq n$. This is the case when $A$ comes from a potential and equation \eqref{1.0} is variational, as \eqref{1.7} below. As a matter of fact, in this note we shall consider only the variational case.}. More precisely, we shall consider \eqref{1.0}  under the condition
\begin{equation}\label{1.2}
g_{1}(x,\snr{z})\mathbb{I}_{n\times n}\le \partial_{z}A(x,z)\le g_{2}(x,\snr{z})\mathbb{I}_{n\times n}
\end{equation}
where the Carath\'eodory-regular functions $g_{1},g_{2}\colon \Omega\times (0,\infty)\to (0,\infty)$ are the smallest and the largest eigenvalues of $\partial_z A$, respectively. To avoid technical complications, unless otherwise stated, we shall assume that $g_{1}$ and $g_{2}$ vanish at most when $\snr{z}=0$, and that $\partial_zA(x, z)$ always exists except when $|z|=0$. The class of equations \eqref{1.0}, subject to \eqref{1.2}, appears in the modelling of a large variety of stationary phenomena such as the theory of electrostatic or electromagnetic potentials, or in the search of vibration modes in elastic structures. We are interested in the regularity of  solutions. Measuring the ellipticity of \eqref{1.0} provides crucial regularity information. To fix ideas, let us consider the linear equation
\begin{equation}\label{1.1}
-\textnormal{div}(\mathcal{A}(x)Du)=0\quad \mbox{in} \ \ \Omega,
\end{equation}
where, for simplicity, $\mathcal{A}$ is a measurable, symmetric and positive definite $\mathbb R^{n\times n}$-matrix. The standard ellipticity quantifier is the ellipticity ratio $\mathcal R$, defined as the ratio of the highest eigenvalue of $\mathcal A$ over the lowest eigenvalue of $\mathcal A$ 
\begin{equation}\label{1.3}
\mathcal{R}(x):=\frac{\mbox{highest eigenvalue of }\mathcal A(x)}{\mbox{lowest eigenvalue of }\mathcal A(x)}\,.
\end{equation}
The condition 
\eqn{lasci}
$$\nr{\mathcal{R}}_{L^{\infty}({\rm B})}\mbox{ is finite for all balls ${\rm B}\Subset \Omega$} $$
  classifies \eqref{1.1} as a uniformly elliptic PDE \cite[Chapter 3]{gt83}; if otherwise $\nr{\mathcal{R}}_{L^{\infty}({\rm B})}$ is not finite for at least one ball ${\rm B}\Subset \Omega$, equation \eqref{1.1} is called nonuniformly elliptic\footnote{Since in this note we are mainly interested in local regularity properties, here we are using suitably localized notions of uniform and nonuniform ellipticity, that slightly differ from the classical ones. Specifically, the finiteness of the ellipticity ratio is analysed only on interior balls, rather than on the whole domain as usually done \cite{gt83}, i.e., $\Omega$ is taken instead of ${\rm B}$ in \eqref{lasci}. This naturally suites to interior regularity. Of course, this distinction is not strictly necessary in our setting.}. In order to  extend this notion to the nonlinear elliptic setting \eqref{1.0}-\eqref{1.2}, let us consider the simplified equation $-\textnormal{div} \ A(Du)=0$ and (formally) differentiate it in the $s$-direction, $s \in \{1, \ldots, n\}$. We find out that every component $v_{s}:=D_{s}u$ solves 
\begin{equation}\label{keepin}
-\textnormal{div}\, (\mathcal{A}(x)Dv_{s})=0,\quad \quad \mathcal{A}(x):=\partial_{z}A(Du(x)),
\end{equation}
i.e., second derivatives of $u$ solve a linear equation with matrix coefficient $\mathcal{A}$ being a nonlinear, measurable function of $Du$. Back to general case, the adaptation of \eqref{1.3} to the nonlinear framework \eqref{1.0} towards gradient regularity is then quite natural \cite{tru67, sim71}. Keeping \eqref{1.2} and \eqref{keepin} in mind, we introduce the pointwise ellipticity ratio for $|z|\not=0$,
\begin{equation}\label{1.4}
\mathcal{R}(x,z):=\frac{g_{2}(x,\snr{z})}{g_{1}(x,\snr{z})},
\end{equation}
and, again for $\snr{z}\not =0$, the nonlocal one \cite{dm21},
\begin{equation}\label{1.5}
\bar{\mathcal{R}}(z,{\rm B}):=\frac{\sup_{x\in {\rm B}}g_{2}(x,\snr{z})}{\inf_{x\in {\rm B}}g_{1}(x,\snr{z})}
\end{equation}
for any ball ${\rm B}\Subset \Omega$. 
Quantities \eqref{1.4}-\eqref{1.5} are well-defined thanks to the ellipticity condition \eqref{1.2} (with the obvious agreement that $\bar{\mathcal{R}}(z,{\rm B})=\infty$ when the  denominator in \eqref{1.5} is zero\footnote{This will never occur in the situations considered here.}). Moreover, by definition, 
\begin{equation}\label{1.6}
\begin{cases}
    \mathcal{R}(x,z)\le \bar{\mathcal{R}}(z,{\rm B})\quad \mbox{for} \  x\in {\rm B}\\
    A(x,z)\equiv A(z) \Longrightarrow   \mathcal{R}(x,z)=\bar{\mathcal{R}}(z,{\rm B})\,.
    \end{cases}
\end{equation}
As suggested in \cite{dm21}, we adopt a slightly more refined taxonomy for uniform ellipticity.  
\begin{itemize}
\item $A$ is uniformly elliptic if
    \begin{equation}\label{unel}
    \sup_{\snr{z}>0}\bar{\mathcal{R}}(z,{\rm B})\lesssim_{{\rm B}} 1,
    \end{equation}
    for all balls ${\rm B}\Subset \Omega$.
    \item $A$ is softly nonuniformly elliptic if it is pointwise uniformly elliptic, i.e., 
    \begin{equation}\label{soft.1}
   \sup_{x\in {\rm B}, \snr{z}>0}\mathcal{R}(x,z)\lesssim_{{\rm B}} 1
    \end{equation}
    holds for all balls ${\rm B}\Subset \Omega$ but
        \begin{equation}\label{soft.2}
    \sup_{\snr{z}>0}\bar{\mathcal{R}}(z,{\rm B})=+\infty
    \end{equation}
    holds for some ball ${\rm B}\Subset \Omega$.
    \item $A$ is (strongly, pointwise) nonuniformly elliptic if
    \begin{equation}\label{strong}
    \sup_{x\in {\rm B}, \snr{z}>0}\mathcal{R}(x,z)=+\infty 
    \end{equation}
    on some ball ${\rm B}\Subset \Omega$\footnote{In line with \cite{sim71}, possible alternative definitions of uniform ellipticity can be given by replacing \eqref{unel} and \eqref{soft.1} with 
    $$
    \limsup_{\snr{z}\to \infty}\bar{\mathcal{R}}(z,{\rm B})<\infty  \quad \mbox{and}\quad    \limsup_{\snr{z}\to \infty} \sup_{x}\mathcal{R}(x,z)<\infty,
    $$ 
respectively. This is according to the fact that one is interested in proving Lipschitz continuity of solutions, after which in several cases the problems in question become uniformly elliptic again. Therefore one is essentially interested in looking at those situations where ellipticity weakens and the gradient becomes simultaneously large.}.
\end{itemize}
\noindent 
The classical definition of nonuniform ellipticity only prescribes to distinguish between the occurrence of \eqref{soft.1} and \eqref{strong} (with ${\rm B}$ replaced by $\Omega$). Instead, here we consider an intermediate notion that on the other hand becomes immaterial in the autonomous case \eqref{1.6}$_2$. Specifically, the pointwise ellipticity ratio captures how the growth of the gradient variable affects ellipticity, while the nonlocal one governs the way space-dependent coefficients mix-up with gradients, and indicates how much coefficients deviate from merely having a perturbative effect. Let us stress that, when dealing with nonuniformly elliptic problems, it is more natural to consider the variational setting, which provides the correct notion of solution (minimizer) in the natural energy setting. Specifically, we consider variational integrals of the type
\begin{equation}\label{fun}
    \mathcal{F}(w,\Omega):=\int_{\Omega}F(x,Dw)\dx,
\end{equation}
governed by an elliptic integrand $F\colon \Omega\times \mathbb{R}^{n}\to [0,+\infty)$. The catch between the functional $\mathcal{F}$ and equation \eqref{1.0} is the Euler-Lagrange equation
\begin{equation}\label{1.7}
-\textnormal{div}\, \partial_{z}F(x,Du)=0\quad \mbox{in} \ \ \Omega, 
\end{equation}
solved by minimizers under reasonable structure conditions. Thanks to \eqref{1.7}, letting $A=\partial_{z}F$, the positions in \eqref{1.4}-\eqref{1.5} can be directly adapted to the variational framework: functions $g_{1}$ and $g_{2}$ in \eqref{1.2} now become  lower and upper bounds for the lowest and the highest eigenvalue of the Hessian $\partial_{zz}F$, respectively. It is convenient to clarify the notion of minimality we are adopting. \begin{definition}[Minima]\label{defi-min} A function $u \in W^{1,1}_{\loc}(\Omega)$ is a \emph{(local) minimizer} of the functional $\mathcal F$ in \eqref{fun} if, for every ball ${\rm B}\Subset \Omega$, $F(\cdot, Du) \in L^1({\rm B})$ and $\mathcal F(u,{\rm B})\leq \mathcal F(w,{\rm B})$ holds for every $w \in u + W^{1,1}_0({\rm B})$. 
\end{definition}
\subsection{Master examples}\label{kex} Basic model examples capturing some of the main aspects of the ellipticity notions discussed above are in the following.
\subsubsection{Uniform ellipticity.} This is the case of equations of the $p$-Laplace type, 
\begin{equation}\label{plap}
\begin{cases}
-\textnormal{div}\left(a(x)\snr{Du}^{p-2}Du\right)=0\quad \mbox{in} \ \ \Omega\\
\displaystyle
1<p<\infty,\quad 1\lesssim a\in L^\infty_{\loc}(\Omega)\,.
\end{cases}
\end{equation}
A direct computation yields\footnote{In this note, with some abuse of notation, we interchange $\sup$ and esssup, the correct meaning being clear from the context.}
\begin{equation*}
\mathcal{R}(x,z)\leq  \frac{\max\{1,p-1\}}{\min\{1,p-1\}} ,\qquad\quad  \bar{\mathcal{R}}(z,{\rm B})\leq\frac{\sup_{x\in {\rm B}}a(x)}{\inf_{x\in {\rm B}}a(x)}  \frac{\max\{1,p-1\}}{\min\{1,p-1\}}, 
\end{equation*}
for all balls ${\rm B}\Subset \Omega$. After the foundational contributions of Ural'tseva \cite{ura68} and Uhlenbeck \cite{uhl77} on the case $a(\cdot)\equiv 1$, maximal regularity for energy solutions\footnote{Energy solutions to \eqref{1.0} with $\snr{A(x,z)}\lesssim \snr{z}^{q-1} +1$, $q>1$, are those distributional solutions $u\in W^{1,1}_{\loc}(\Omega)$ satisfying the additional requirement $u \in W^{1,q}_{\loc}(\Omega)$, while distributional solutions are only required to satisfy $A(\cdot, Du)\in L^1_{\loc}(\Omega;\mathbb R^n)$. These are sometimes called very weak solutions and might be irregular \cite{ser64}. In the case \eqref{plap} an energy solution $u$ is distributional solution such that  $u \in W^{1,p}_{\loc}(\Omega)$.} to \eqref{plap}, i.e., gradient H\"older continuity, was established in the works of   DiBenedetto \cite{dib83} and Manfredi \cite{man88}, see also Kuusi \& Mingione \cite{km14} for related borderline cases.
\subsubsection{Soft nonuniform ellipticity.} A significant example is the double phase integral
\begin{equation}\label{dp}
\begin{cases} 
    \displaystyle
   \,  \mathcal{D}(w,\Omega):=\int_{\Omega}\snr{Dw}^{p}+\aaa(x)\snr{Dw}^{q}\dx\\[10pt]\displaystyle
  \,   1<p< q<\infty,\quad 0\le \aaa\in L^{\infty}_{\loc}(\Omega).
\end{cases}
\end{equation}
In this case we have
 \begin{equation}\label{dp.1}
\begin{cases}
\,\displaystyle  \mathcal{R}(x,z)\leq \frac{q\max\{1,q-1\}}{p\min\{1,p-1\}} \\ \\
\,  \displaystyle \bar{\mathcal{R}}(z,{\rm B})\leq  \frac{q\max\{1,q-1\}}{p\min\{1,p-1\}} \left(\nr{\aaa}_{L^{\infty}({\rm B})}\snr{z}^{q-p}+1\right).
  \end{cases}
\end{equation}
By \eqref{dp.1}$_1$ the integrand governing the functional $\mathcal{D}$ is uniformly elliptic in the classical sense, i.e. \eqref{soft.1} holds. As a matter of fact, considering $\aaa(x)=|x|^\alpha$ we have $\bar{\mathcal{R}}(z,\mathcal B_r)\approx r^\alpha \snr{z}^{q-p}+1$ so that \eqref{soft.2} holds.
The possible blow-up of $\bar{\mathcal{R}}$ indicates that the vanishing of coefficient $\aaa$ is the sole responsible of the (mild) nonuniformity in \eqref{dp}. The functional $\mathcal{D}$ was introduced by Zhikov \cite{zhi87, jko94} in the setting of homogenization, and to investigate the occurrence of Lavrentiev phenomenon. The mix-up of the coefficient $\aaa$ and gradient variable in \eqref{dp} might lead to the formation of singularities even in scalar, nondegenerate problems, as shown by Esposito \& Leonetti \& Mingione \cite{elm04}, Fonseca \& Mal\'y \& Mingione \cite{fmm04}, and Balci \& Diening \& Surnachev \cite{bds20,bds23}. A complete regularity theory was eventually obtained by Baroni \& Colombo \& Mingione \cite{cm15,bcm18}. We will further discuss \eqref{dp} in Section \ref{1.xx}. Another popular example of softly nonunformly elliptic functional is the variable exponent one
\eqn{px}
$$
w \mapsto \int_{\Omega} |Dw|^{\mf{p}(x)}\dx
$$
\noindent where $\mf{p}\colon \Omega \mapsto (1, \infty)$  is a measurable function such that $1< p \leq \mf{p}(\cdot)\leq q <\infty$, and for which we have  
$$
\begin{cases}
\, \displaystyle  \mathcal R (x,z)\leq \frac{q\max\{1,q-1\}}{p\min\{1,p-1\}}\\ \\
\, \displaystyle \bar{\mathcal{R}} (z,{\rm B}) \leq \frac{q\max\{1,q-1\}}{p\min\{1,p-1\}}\max\{ \snr{z}, 1/\snr{z}\}^{\textnormal{osc}_{{\rm B}}\mf{p}}
\end{cases}
$$
Maximal regularity for minima of functionals modelled by the one in \eqref{px} was established by Acerbi \& Mingione in \cite{am01}. Specifically, for minimizers 
$u$
of \eqref{px}, a Schauder-type result holds, asserting that the gradient $Du$ is locally Hölder continuous whenever the exponent function $\mf{p}$ itself is Hölder continuous.
More examples of softly nonuniformly elliptic functionals are discussed in \cite{ho22a,bb25}.\subsubsection{Strong nonuniform ellipticity.} A mostly celebrated example in this class is given by the area integral\footnote{Beware! In this case the integrand in question has linear growth so that the natural ambient function space is BV  and the problems must be formulated accordingly, as for instance in \cite{gms79a, gms79b}. Therefore the discussion here is formal, and only regards the ellipticity ratio. We refer to \cite{bs13, bs15} for larger discussion on functionals with linear growth.},
 \begin{equation}\label{ms}
   w \mapsto \int_{\Omega}\sqrt{1+\snr{Dw}^{2}}\dx,
 \end{equation}
 and the related Euler-Lagrange equation, i.e., the minimal surface equation
 $$
 -\textnormal{div} \left(\frac{Du}{\sqrt{1+\snr{Du}^2}}\right)=0
 $$
intensively studied for instance in the seminal work of Bombieri \& De Giorgi \& Miranda \cite{bdm69}, Ladyzhenskaya \& Ural'tseva \cite{lu70}, Simon \cite{sim76}, Giusti \cite{giu78}, Giaquinta \& Modica \& Souček \cite{gms79b}, and Trudinger \cite{tru81}, see also Bildhauer \& Fuchs \cite{bf03}, Beck \& Schmidt \cite{bs13,bs15} and Gmeineder \& Kristensen \cite{gk19,gk24} for recent advances and more general structures. The pointwise ellipticity ratio of \eqref{ms} (that coincides with the nonlocal one being the functional autonomous)  features a quadratic rate of blow-up, that is $\mathcal{R}(z)= 1+\snr{z}^{2}$. Analogous considerations hold for more general area type functionals:
$$
   w \mapsto\int_{\Omega}(1+\snr{Dw}^{m})^{1/m}\dx,\qquad m>1,
$$
for which  $\mathcal{R}(z)\approx 1+\snr{z}^{m}$ for sufficiently large $|z|$, see \cite{gms79b,bs15}. Variational integrals at linear growth are far from being the only strongly nonuniformly elliptic examples. Strong rates of nonuniformity are indeed typical of several functionals arising through various fields, such as fluid dynamics, materials science or nonlinear elasticity. Examples are given by slow-growing functionals \cite{fs99,fm00,bf03} (nearly linear growth conditions),
\begin{equation}\label{llog}
w \mapsto \int_{\Omega}\snr{Dw}\log(1+\snr{Dw})\dx,
\end{equation}
satisfying
$
    \mathcal{R}(z)\approx \log(1+\snr{z})
$
for $\snr{z}$ large, 
or convex polynomials,
\begin{equation}\label{cp}
\begin{cases}
\displaystyle
\, w \mapsto \int_{\Omega}\snr{Dw}^{p}+\sum_{i=1}^{n}\aaa_{i}(x)\snr{D_iw}^{q_{i}}\dx\\[10pt]\displaystyle
\, 1<p\leq q_{1}\leq \cdots\leq q_{n}<\infty\\[12pt]\displaystyle
\, 0\leq \aaa_{i}(\cdot)\in L^{\infty}_{\loc}(\Omega), \quad  i\in \{1,\cdots,n\},
\end{cases} 
\end{equation}
in which case we find 
\begin{equation*}
\mathcal{R}(x,z)+ \bar{\mathcal{R}}(z,{\rm B})\lesssim_{p,q_{i},\nr{\aaa_{i}}_{L^{\infty}({\rm B})}}\snr{z}^{q_{n}-p} +1,
\end{equation*}
for $x \in {\rm B}$, see \cite{mar89,elm04,fgk04,bb20,dkk24}; or integrands at fast exponential growth, i.e.,
\begin{equation*}
\begin{cases}
\displaystyle
  \,   w\mapsto \int_{\Omega}\exp\left(a(x)\snr{Dw}^{p}\right)\dx\\[10pt]\displaystyle
  \,   1<p<\infty,\quad 1\lesssim a\in L^{\infty}_{\loc}(\Omega),
    \end{cases}
\end{equation*}
where
\begin{equation*}
\mathcal{R}(x,z)\approx \snr{z}^{p}+1,\qquad\quad   \bar{\mathcal{R}}(z,{\rm B})\approx   \left(\snr{z}^{p}+1\right)e^{(\textnormal{osc}_{{\rm B}}a)\snr{z}^{p}},
\end{equation*}
for $\snr{z}$ large; see \cite{mar96,bm20,dm21}. A novel class of integrals, first considered in \cite{dm23b}, and appearing as a combination of classical nearly linear growth ones \eqref{llog} and weighted power terms, is that of log-double phase problems, i.e.,
\eqn{limiting}
$$
\begin{cases}
\displaystyle
\, w\mapsto \int_{\Omega}\snr{Dw}\log(1+\snr{Dw})+\aaa(x)\snr{Dw}^{q}\dx\\[10pt]\displaystyle
\, 1<q<\infty,\quad 0\le \aaa\in L^{\infty}_{\loc}(\Omega).
\end{cases}
$$ 
The related ellipticity ratios behave as hybrids of those of \eqref{dp} and \eqref{llog}, i.e.,
\begin{equation*}
\begin{cases}
\, \mathcal{R}(x,z)\lesssim_{q} \log(1+\snr{z}) +1\\
\, \bar{\mathcal{R}}(z,{\rm B})\lesssim_{q} \log(1+\snr{z})+\nr{\aaa}_{L^{\infty}({\rm B})}\snr{z}^{q-1}+1.
\end{cases}
\end{equation*}
The functional in \eqref{limiting}, readable as a limiting configuration of \eqref{dp} as $p\to 1$ while approaching nearly linear growth conditions as in \eqref{llog}, leads to sharp outcomes on the validity of Schauder theory for anisotropic problems whose growth is arbitrarily close to linear \cite{dm23b,ddp24,dp24}, and offers deep insights on subtle singularity phenomena for more general nonautonomous area-type integrals. 
\section{Polynomial nonuniform ellipticity}\label{pne} A mostly common rate of nonuniformity for equations and functionals is the polynomial one, meaning that the pointwise ellipticity ratio behaves (at infinity) as a positive power of the gradient variable, i.e.,
\begin{equation}\label{2.0}
\mathcal{R}(x,z)\approx_{\snr{z}\ge 1} \snr{z}^{\delta},\qquad \quad \delta>0.
\end{equation}
See for instance the foundational regularity works of Serrin \cite{ser69}, Ladyzhenskaya \& Ural'tseva \cite{lu70}, Ivanov \cite{iva72}, Simon \cite{sim76}, Trudinger \cite{tru67, tru81} and Ural'tseva \& Urdaletova \cite{uu84}. Note that many of such results rely on the existence of strong solutions and on the availability of total derivatives, the latter requiring smoothness assumptions on coefficients that automatically rule out the classical Schauder setting. Marcellini's innovative variational approach \cite{mar86,mar89,mar91} to nonuniformly elliptic equations \eqref{1.0} marked significant progress. The key condition in Marcellini's work is trapping the lowest and highest eigenvalues of $\partial_{z}A$ between two different powers of the gradient variable in a way to extend the $p$-Laplace type behaviour, i.e., 
\begin{equation}\label{pq}
\begin{cases}
\displaystyle
\, \snr{z}^{p-2}\mathbb{I}_{n\times n}\lesssim \partial_{z}A(x,z)\lesssim \snr{z}^{q-2}\mathbb{I}_{n\times n} \ \ \mbox{in the sense of matrices,}\\[10pt]\displaystyle
\, \mbox{for some} \ \ 1<p\le q \ \ \mbox{and all} \ \ \snr{z}\ge 1,\ x\in \Omega,
\end{cases}
\end{equation}
the growth at infinity of the related ellipticity ratio can be controlled in terms of the difference $q-p$,
\begin{equation}\label{pqpq}
\mathcal{R}(x,z)\lesssim\snr{z}^{q-p}\,,  \quad |z|\geq 1.
\end{equation}
This naturally relates to \eqref{2.0}. 
The growth rate of $\mathcal{R}$ can be tamed by taking $q-p$ sufficiently small, so that  the asymptotic in \eqref{pqpq} suggests that only a moderate blow-up rate of $\mathcal{R}$ gives hope for regular solutions.\footnote{Conditions \eqref{pq}-\eqref{pqpq} can be recasted in the variational setting by replacing $A$ with $\partial_{z}F$, and $\partial_{z}A$ with $\partial_{zz}F$.} In fact, building on earlier works of Marcellini \cite{mar89,mar91} and Giaquinta \cite{gia87}, Min-Chun \cite{mc92} proved that for $n\ge 6$, the function
\begin{equation}\label{mc1}
    u(x):=\sqrt{\frac{n-4}{24}}\frac{x_{n}^{2}}{\sqrt{\sum_{i=1}^{n-1}x_{i}^{2}}}- \frac{2}{n-2} \sqrt{\frac{n-4}{24}\sum_{i=1}^{n-1}x_{i}^{2}},
\end{equation} 
which is unbounded on the  line $(0,\cdots,0,x_{n})$, minimizes the convex polynomial
\begin{equation}\label{mc2}
w \mapsto \int_{\Omega}\snr{Dw}^{2}+\frac{1}{2}\snr{D_{n}w}^{4}\dx.
\end{equation}
The functional in \eqref{mc2}, which is of the type in \eqref{cp}, highlights a peculiarity of nonuniformly elliptic problems: already for strongly convex, autonomous functionals 
\begin{equation}\label{fun0}
\begin{cases}
\displaystyle
\, \mathcal{F}_{0}(w,\Omega):=\int_{\Omega}F_{0}(Dw)\dx\\[8pt]\displaystyle
\, \snr{z}^{p}\lesssim F_{0}(z)\lesssim 1+\snr{z}^{q},\quad \quad  1<p\le q<\infty,
\end{cases}
\end{equation}
satisfying \eqref{pq} with $\partial_{z}A\equiv \partial_{zz}F_{0}$, minimizers may be unbounded when a constraint of the type
\begin{equation}\label{pq0}
\frac{q}{p}<1+ \texttt{o}_n,\qquad \quad  \lim_{n} \texttt{o}_n=0
\end{equation}
is violated. This is an element of sharp contrast with the uniformly elliptic setting. In fact, uniform ellipticity allows the formation of singularities only in vectorial problems \cite{deg68,maz68,nec77,sy02,ms16}, while Marcellini, Giaquinta and Min-Chun's counterexamples are scalar. In particular, the regularity of minima of \eqref{fun0} does not hold for arbitrary choices of $p$ and $q$, that actually must satisfy a precise quantitative constraint involving exponents $(p,q)$ and the ambient dimension $n$. The optimal bound in \eqref{pq0} is known only at the 0-order scale: Marcellini \cite{mar91} exhibited an unbounded minimizer of an integral of type \eqref{fun0} provided that
\begin{equation}\label{2.1}
q>\frac{p(n-1)}{n-1-p}\qquad \mbox{and}\qquad 1<p<n-1
\end{equation}
holds, which is consistent with the construction in \eqref{mc1}-\eqref{mc2}: dimensional condition $n\ge 6$ is indeed equivalent to $4>2(n-1)/(n-3)$. Local boundedness of minima  in violation of \eqref{2.1}, i.e.,
\begin{equation*}
\left\{
\begin{array}{c}
\displaystyle
1<q\le \frac{p(n-1)}{n-1-p}\qquad \mbox{if} \ \ 1<p<n-1\\[12pt]\displaystyle
1<p\le q<\infty\qquad \mbox{if} \ \ p\ge n-1
\end{array}
\right.
\end{equation*}
was obtained thirty years later by Hirsch \& Schäffner \cite{hs21} via fine optimization techniques, also relying on methods by Bella \& Schäffner \cite{bs20}. Nonetheless, the analogous issue at the gradient level, in full generality, is still open. Marcellini \cite{mar91} proved gradient higher differentiability and local Lipschitz continuity for minima of \eqref{fun0} with \eqref{pq} in force\footnote{See \cite{mar91} for precise assumptions, or see \eqref{assf} below.} if 
\begin{equation*}
 \frac{q}{p}<1+\frac{2}{n},   
\end{equation*}
a constraint later on updated to
\begin{equation*}
\frac{q}{p}<1+\frac{2}{n-1}
\end{equation*}
by Bella \& Schäffner \cite{bs20,bs24,sch24}. In case $F_{0}$ is an elliptic convex polynomial, the above restrictions were improved to
\begin{equation}\label{2.3}
\left\{
\begin{array}{c}
\displaystyle
2\le p\le q\le \frac{p(n-1)}{n-3}\qquad \mbox{if} \ \ n\ge 4\\[10pt]\displaystyle
2\le p\le q<\infty\qquad \mbox{if} \ \ n\in \{2,3\}
\end{array}
\right.
\end{equation}
by Koch \& Kristensen and the author \cite{dkk24}. When $p=2$, \eqref{2.3} matches the threshold determined in \cite{gia87,mar89,mc92} and yields sharp gradient regularity; the general case $p\not=2$ is still far from being understood. Aside from the optimal bound under which gradient regularity can be achieved, the theory for autonomous, nonuniformly elliptic functionals is at a very advanced stage. It is therefore natural to wonder what happens when external ingredients come into play, specifically if Schauder theory still holds for nonuniformly elliptic problems. For a further discussion of nonuniformly elliptic problems and techniques we refer to \cite{def25, dm25b}. 
\section{Uniformly elliptic Schauder theory}\label{ues}
Consider the linear equation with variable coefficients \eqref{1.1}, 
where the matrix $  \mathcal A$ is symmetric (for simplicity), bounded and elliptic, i.e., $\mathbb{I}_{n\times n}\leq \mathcal A\in L^\infty$. Schauder theory claims
\begin{equation}\label{3.0} 
    \mathcal A\in C^{0,\alpha}_{\loc}(\Omega;\mathbb{R}^{n\times n}) \ \Longrightarrow \ Du\in C^{0,\alpha}_{\loc}(\Omega;\mathbb{R}^{n}),\qquad \alpha\in (0,1).
\end{equation}
That is, the gradient of the solution inherits the same degree of regularity of coefficients; this is optimal (and fails in the borderline case $\alpha=1$).  
This kind of results were  pioneered by Hopf \cite{hop28}, Giraud \cite{gir29}, Caccioppoli \cite{cac34}, and Schauder \cite{sch34a,sch34b}, in various forms and today are known as Schauder estimates, with parabolic,  fully nonlinear and vectorial analogs \cite{adn59, adn64,cc95}. Schauder theory is a basic tool in elliptic and parabolic PDEs and in the Calculus of Variations, influencing diverse areas of analysis such as nonlinear diffusion, potential theory, field theory or differential geometry. The classical outcome is twofold: 
\begin{itemize}
    \item regularity estimates of type \eqref{3.0}, local or global, establishing the principle that solutions to nonautonomous equations are as regular as the ingredients allow;
    \item boundedness or compactness of the inverse of certain elliptic operators;
\end{itemize}
see \cite{kic06}. The latter allows proving existence of solutions via fixed point theorems, that was a standard  argument before the introduction of energy methods and Sobolev spaces. Here, we shall be interested in the former. The results in \cite{hop28,gir29,cac34,sch34a,sch34b} rely on potential theoretic techniques. However, more direct approaches are today available avoiding the use of representation formulas, as we shall see in a few lines directly in the nonlinear case. For instance,  Campanato \cite{cam65} used  suitable function spaces; Trudinger \cite{tru86} used convolution techniques, and Leon Simon \cite{sim97} used blow-up methods. The nonlinear theory dates back to the 1980s, with the work of Giaquinta \& Giusti \cite{gg82,gg83,gg84}, DiBenedetto \cite{dib83} and Manfredi \cite{man88}. The underlying principle is essentially to use the fact that H\"older continuity of coefficients allows to prove the closeness of the original solution $u$ to to that of "frozen" problems of the type
$$
    \begin{cases}
 \,   -\textnormal{div}(\mathcal A(x_{0})Dv)=0\quad &\mbox{in} \ \ B_{r}(x_{0})\\
 \,   \ v=u\quad &\mbox{on} \ \ \partial B_{r}(x_{0}),
    \end{cases}
$$
on small balls $B_{r}(x_{0})\Subset \Omega$. Let's see an application of this principle directly in the nonlinear case
\eqn{nona0}
$$
 -\textnormal{div}((\mathcal A(x)Du\cdot Du)^{(p-2)/2}\mathcal A(x)Du)=0,
$$
where we consider an energy solution $u\in W^{1,p}_{\loc}(\Omega)$, $p>1$; note that \eqref{nona0} reduces to \eqref{1.1} for $p=2$. For a given ball $B_{r}(x_{0})\Subset \Omega$, the related nonlinear lifting (balayage) can be defined as the solution to Dirichlet problem
$$
   \begin{cases}
  \,   -\textnormal{div}((\mathcal A(x_0)Dv\cdot Dv)^{(p-2)/2}\mathcal A(x_0)Dv)=0\quad &\mbox{in} \ \ B_{r}(x_{0})\\
   \,   v=u\quad &\mbox{on} \ \ \partial B_{r}(x_{0}).
   \end{cases} 
$$
Thanks to Ural'tseva-Uhlenbeck theory \cite{ura68,uhl77}, $Dv$ is locally H\"older continuous in $B_{r}(x_{0})$, and 
\eqn{3.5}
$$
\begin{cases}
\displaystyle \nr{Dv}_{L^{\infty}(B_{r/2}(x_0))}^p \lesssim_{\nr{\mathcal A}_{L^{\infty}(B_{r}(x_{0}))}} \mint_{B_{r}(x_{0})}\snr{Dv}^{p}\dx\\
\displaystyle \mint_{B_{\sigma}(x_{0})}\snr{Dv-(Dv)_{B_{\sigma}(x_{0})}}^{p}\dx\\
 \displaystyle\qquad  \lesssim_{n,p,\nr{\mathcal A}_{L^{\infty}(B_{r}(x_{0}))}} \left(\frac{\sigma}{\varrho}\right)^{\beta_{0}p}\mint_{B_{\varrho}(x_{0})}\snr{Dv-(Dv)_{B_{\varrho}(x_{0})}}^{p}\dx
 \end{cases}
$$
hold on all concentric balls $B_{\sigma}(x_{0})\subset B_{\varrho}(x_{0})\Subset B_{r}(x_{0})$ (see for instance \cite{man88}, \cite[Theorem 3.2]{km18}). Here $\beta_0\equiv \beta_0(n,p)\in (0,1)$ is a universal H\"older continuity exponent, which can be taken to be $\beta_0=1$\footnote{In the linear case this is in fact a classical result of Campanato \cite{cam65} and \eqref{3.5}$_2$ holds with $\beta_0=1$. When $p\not=2$ estimate \eqref{3.5}$_2$ is a rigid form of Ural'tseva-Uhlenbeck theory, and one cannot take $\beta_0=1$ in general,  by counterexamples.} when $p=2$. Moreover, the $\alpha$-H\"older continuity of the matrix  $\mathcal A$ grants the validity of comparison estimate
\begin{equation}\label{3.6}
\mint_{B_{r}(x_{0})}\snr{Du-Dv}^{p}\dx\lesssim_{n,p,[\mathcal A]_{0,\alpha;B_{r}(x_{0})}}r^{\alpha\min\{2,p\}}\mint_{B_{r}(x_{0})}\snr{Du}^{p}\dx.
\end{equation}
The combination of \eqref{3.5}-\eqref{3.6} allows transferring smoothness from $v$ to $u$
\eqn{nona}
$$
\mint_{B_{r}(x_0)}\snr{Du-(Du)_{B_{r}}}^{p}\dx \lesssim r^{p\beta} , \qquad \beta:=\min\{\beta_0, \alpha\min\{2/p,1\}\}
$$
for all sufficiently small $\varrho$, that implies that 
\eqn{brings}
$$Du \in C^{0, \beta}_{\loc}(\Omega;\mathbb R^n)$$ via a classical integral characterization due to Campanato and Meyers. Note that in the linear case $p=2$, we can take $\beta_0=1$ and therefore   $\beta=\alpha$. This allows to recover \eqref{3.0}. The very same circle of ideas applies to equations as in \eqref{1.0} satisfying the same growth/ellipticity and oscillation features as the models in \eqref{plap} and \eqref{nona0}
, and to nonautonomous problems governed by general uniformly elliptic operators  \cite{lie91}. 
\begin{remark}\label{rem1}
\emph{We identify three aspects of the uniformly elliptic Schauder theory developed between the early 1930s and the end of the 1980s. 
\begin{itemize}
    \item Proof by perturbation. That is, use of the H\"older continuity of coefficients to be quantitatively close to a uniformly elliptic problem with constant (frozen) coefficients \eqref{3.6}, whose solutions (liftings) enjoy homogeneous excess decay estimates \eqref{3.5}, and then exploit homogeneity to recover similar integral decay estimates for solutions to the original nonautonomous problem \eqref{nona}.
    \item Smooth data never obstruct the regularity of energy solutions \eqref{3.0}.
    \item Lipschitz bounds via gradient H\"older continuity. In the classical theory of uniformly elliptic equations the standard way to achieve gradient $L^\infty$-bounds in nonautonomous, uniformly elliptic problems with H\"older coefficients as \eqref{nona0} is to prove first gradient H\"older continuity for some exponent, and then retain local boundedness. 
\end{itemize}}
\end{remark}
The last point of Remark \ref{rem1} is of particular interest here. Indeed, when the (nonlinear) equation in question is nondegenerate or nonsingular, the Lipschitz information is employed afterwards to in a sense  linearize the equation and establish that the gradient of solutions shares the same H\"older exponent as the coefficient. As an example, consider the nondegenerate $p$-Laplace equation
\begin{equation}\label{plapnd}
\begin{cases}
\,  -\textnormal{div}(a(x)(\snr{Du}^2+1)^{(p-2)/2} Du)=0\quad \mbox{in} \ \ \Omega\\[10pt]\, \displaystyle \, 
1<p<\infty,\quad 1\lesssim a\in C^{0,\alpha}_{\loc}(\Omega), \quad \alpha \in (0,1)\,.
\end{cases}
\end{equation}
In this case the technique outlined in this section leads to \eqref{brings}. Further regularity techniques, explained in \cite[Section 5]{dm25a}, eventually lead to 
\eqn{bring2}
$$Du \in C^{0, \beta}_{\loc}(\Omega;\mathbb R^n)\Longrightarrow Du \in L^{\infty}_{\loc}(\Omega;\mathbb R^n)\Longrightarrow Du \in C^{0, \alpha}_{\loc}(\Omega;\mathbb R^n)$$
and therefore to the same sharp result of the linear case notwithstanding $p\not=2$\footnote{In fact, already in the degenerate case \eqref{nona0}, the standard derivation of \eqref{brings} goes as follows: in a first step one combines \eqref{3.5}$_1$ and \eqref{3.6} to obtain that $Du$ belongs to a suitable Morrey space. Then ones uses this to combine this time \eqref{3.5}$_2$ and \eqref{3.6} to get that $Du$ is locally H\"older continuous with some exponent \cite{man88}. In particular, this implies that $Du$ is locally bounded. Finally, one use this last information to recombine in a different way \eqref{3.5}$_2$ and \eqref{3.6} to get \eqref{brings}. The derivation of \eqref{bring2} in the non-degenerate case \eqref{plapnd} is more delicate but follows similar steps \cite[Section 5]{dm25a}.}. Note that the bootstrap in \eqref{bring2} only works in the nondegenerate/singular case \eqref{plapnd}. We also stress the fact that the second implication in \eqref{bring2} works also when the equations considered are nonuniformly elliptic (but still non degenerate), as, when the gradient is known to be bounded, equations in a sense become uniformly elliptic again. As remarked above, the first implication in \eqref{bring2} is instead essentially the only way to get gradient $L^\infty$-bounds in presence of H\"older coefficients in the classical theory. We will turn back on these points in Section \ref{2.xx}, where in fact this way is reversed.
\section{Nonuniformly elliptic Schauder theory}
\noindent Here the focus will be on nonautonomous variational integrals of the type \eqref{fun} verifying 
\begin{equation}\label{assf}
    \begin{cases}
        \, \bar{z}\mapsto F(x,\bar{z})\in C^{2}(\mathbb{R}^{n}\setminus \{0\})\cap C^{1}(\mathbb{R}^{n})\vspace{0.8mm}\\
        \, \snr{z}^{p}\lesssim F(x,z)\lesssim 1+\snr{z}^{q}\vspace{0.8mm}\\
        \, \langle\partial_{zz}F(x,z)\xi,\xi\rangle\gtrsim (\mu^{2}+\snr{z}^{2})^{\frac{p-2}{2}}\snr{\xi}^{2}\vspace{0.8mm}\\
        \, \snr{\partial_{zz}F(x,z)}\lesssim (\mu^{2}+\snr{z}^{2})^{\frac{p-2}{2}}+(\mu^{2}+\snr{z}^{2})^{\frac{q-2}{2}}\vspace{0.8mm}\\
        \, \snr{\partial_{z}F(x_{1},z)-\partial_{z}F(x_{2},z)}\lesssim \snr{x_{1}-x_{2}}^{\alpha}\left(1+\snr{z}^{q-1}\right),
    \end{cases}
\end{equation}
for all $z\in\mathbb{R}^{n}\setminus \{0\}$, $\xi\in \mathbb{R}^{n}$, $x,x_{1},x_{2}\in \Omega$, and some $\mu\in [0,1]$\footnote{The parameter $\mu$ serve to distinguish the degenerate case $\mu=0$ from the non-degenerate one $\mu>0$. We refer to \cite{dm25a} for comments on the assumptions \eqref{assf}.}, $\alpha\in(0,1]$; as usual $1< p \leq q $. Note that under such assumptions Definition \ref{defi-min} implies that any minimizer automatically belongs to $W^{1,p}_{\loc}(\Omega)$. Moreover, \eqref{assf} imply the following growth bound on the (pointwise) ellipticity ratio:
$$
\mathcal{R}(x,z)\lesssim \snr{z}^{q-p}+1\,.
$$
Note that functionals \eqref{dp}, \eqref{px}, \eqref{cp}, \eqref{sm}, are all specific instances of \eqref{assf} for suitable choices of $p,q,\mu$ and H\"older continuous dependence on $x$.   
For nonautonomous integrals as \eqref{fun} subject to \eqref{assf} and $p<q$, the regularity of minima unavoidably finds an obstruction in the possible occurrence of Lavrentiev phenomenon, that is, for instance, 
\begin{equation}\label{5.1} 
\inf_{w\in u_{0}+W^{1,p}_{0}(B)}\mathcal{F}(w,B)<\inf_{w\in u_{0}+W^{1,q}_{0}(B)}\mathcal{F}(w,B)
\end{equation}
with $u_{0}\in W^{1,\infty}(B)$ for some ball $B\Subset \Omega$\footnote{By standard density arguments \eqref{5.1} does not occur when $p=q$.}. In other words, as shown in \cite[Section 3]{elm04}, in connection to \eqref{5.1} it is possible that minima $u$ fail even the basic regularity upgrade $u\in W^{1,q}_{\loc}(\Omega)$. The occurrence of phenomena like \eqref{5.1} leads to consider a larger class of functionals based on lower semicontinuous envelopes built along sequences of more regular functions. The Lebesgue-Serrin-Marcellini (LSM)
relaxation (extension) $\bar{\mathcal{F}}$ of functional $\mathcal{F}$ in \eqref{fun} is defined as 
\begin{equation}\label{LSM}
\bar{\mathcal{F}}(w,U):=\inf_{\{w_{i}\}\subset W^{1,q}(U)}\left\{\liminf_{i\to \infty}\mathcal{F}(w_{i},U)\colon w_{i}\rightharpoonup w \ \ \mbox{in} \ \ W^{1,p}(U)\right\}
\end{equation}
whenever $U\subset \Omega$ is an open subset.  Quantifying the distance between $\mathcal{F}$ and $\bar{\mathcal{F}}$ leads to define the Lavrentiev gap functional $\mathcal{L}_{\mathcal{F}}$ \cite{bm92} as
$$
\mathcal{L}_{\mathcal{F}}(w,B):=\begin{cases}
\ \bar{\mathcal{F}}(w,B)-\mathcal{F}(w,B)\quad &\mbox{if} \ \ \mathcal{F}(w,B)<\infty\vspace{0.8mm}\\
\ 0\quad &\mbox{if} \ \ \mathcal{F}(w,B)=\infty,
\end{cases}
$$
for every ball $B\Subset \Omega$, 
which in a sense provides a quantitative measure of the occurrence of phenomena like \eqref{5.1}. Note that
\begin{itemize}
\item $\mathcal{F}(\cdot, U)\le \bar{\mathcal{F}}(\cdot, U)$. Indeed, the convexity of $z\mapsto F(\cdot,z)$ implied by $\eqref{assf}_{3}$ grants the $W^{1,p}$-weak lower semicontinuity of $\bar{\mathcal{F}}$. It follows that $\mathcal{L}_{\mathcal{F}}$ is non-negative.  
\item $\mathcal{F}(\cdot, B)= \bar{\mathcal{F}}(\cdot, B)$ for every ball $B\Subset \Omega$, when $\mathcal F$ is autonomous, i.e., the integrand $F$ does not depend on $x$.  This follows from a simple convolution argument \cite[Lemma 12]{elm04}. 
\item  $\mathcal F(\cdot, U)=\bar{\mathcal F}(\cdot, U)$ on $W^{1,q}(U)$. Therefore $\bar{\mathcal F}$ can be thought as an extension of $\mathcal F$, when this last functional is initially defined on the smaller space $W^{1,q}(U)$, to the whole $W^{1,p}(U)$. This explains the use of the word extension when referring to $\mathcal F$\footnote{The LSM relaxation was first introduced by Marcellini in \cite{mar86,mar89b} to describe cavitation phenomena in Nonlinear Elasticity, and it is connected to similar, earlier constructions of Lebesgue \cite{leb02} and Serrin \cite{ser61}. It was eventually intensively developed in the literature \cite{mar89,fm97a,fm97b,bfm98,elm04,sch09,def22,dm23a,ds23,gk24,ddp24}. Variants of the definition in \eqref{LSM} are possible, and, often, useful. For instance, the following local version is largely studied in the literature:
$$
\bar{\mathcal{F}}_{\loc}(w,U):=\inf_{\{w_{i}\}\subset W^{1,q}_{\loc}(U)}\left\{\liminf_{i\to \infty}\mathcal{F}(w_{i},U)\colon w_{i}\rightharpoonup w \ \ \mbox{in} \ \ W^{1,p}(U)\right\}.
$$
For both $\bar{\mathcal{F}}$ and $\bar{\mathcal{F}}_{\loc}$ very interesting measure representation results are available \cite{fm97a}. We refer to \cite{sch09} for further variants and related results. Most of the results for $\bar{\mathcal{F}}$ presented here extend to $\bar{\mathcal{F}}_{\loc}$.}.
\end{itemize}
The LSM relaxation naturally comes along with a notion of minimality. 
\begin{definition}[Relaxed Minima]\label{defi-min2} A function $u \in W^{1,p}_{\loc}(\Omega)$ is a \emph{(local) minimizer} of $\bar{\mathcal F}$ in \eqref{LSM} if,  for every ball ${\rm B}\Subset \Omega$, $\bar{\mathcal F}(u,{\rm B})<\infty $ and $\bar{\mathcal F}(u,{\rm B})\leq \bar{\mathcal F}(w,{\rm B})$ holds for every $w \in u + W^{1,p}_0({\rm B})$. 
\end{definition}
The approach used in \cite{elm04}, eventually widely adopted afterwards in the literature, prescribes to prove regularity of (local) minimizers of $\bar{\mathcal F}$\footnote{This, in a sense, automatically excludes the occurrence of phenomena like \eqref{5.1} and allows to concentrate on a priori estimates for more regular solutions.} rather than those of $\mathcal F$.  But then 
\begin{itemize}
\item When 
\eqn{4.2}
$$\frac qp<1+\frac \alpha n $$ minima of $\bar{\mathcal F}$ also (locally) minimize $\mathcal F$. Therefore $\bar{\mathcal F}$ can be used to detect regular minima of $\mathcal F$. In a sense, 
drawing a parallel with the classical theory of elliptic equations, minima of $\mathcal F$ detected by $\bar{\mathcal F}$ correspond to distributional solutions to \eqref{1.0} that are also energy solutions.  See Theorem \ref{t1}.
\item When the Lavrentiev gap functional vanishes on a  minimizer $u$ of the original functional $\mathcal F$, that is, when $\mathcal{L}_{\mathcal F}(u,{\rm B})=0$ holds for all balls ${\rm B}\Subset \Omega$, then 
\eqn{newminima}
$$
\bar{\mathcal{F}}(u, {\rm B}) = \mathcal{F}(u,{\rm B})+\mathcal{L}_{\mathcal F}(u,{\rm B}) = \mathcal{F}(u,{\rm B}) \leq \mathcal{F}(w,{\rm B})\leq \bar{\mathcal{F}}(w, {\rm B}) 
$$
holds whenever $w\in u+W^{1,p}_0({\rm B})$. Therefore $u$ is also a minimizer of the $\bar{\mathcal F}$ and its regularity properties can be derived from the regularity results available for minima of $\bar{\mathcal F}$.  See Corollary \ref{t1.2}. The Lavrentiev gap actually vanishes in a variety of situations, as for instance detailed in \cite[Section 5]{elm04}. See Corollary \ref{main}.
  \end{itemize}
Within this framework, specifically thought for nonautonomous functionals, we shall describe key aspects in a rich regularity  theory initiated by the seminal work of Esposito \& Leonetti \&  Mingione \cite{elm04}. This theory explores the interplay between coefficients regularity, the growth behavior of the ellipticity ratio, and the presence of the Lavrentiev phenomenon, ultimately culminating in the recent establishment of Schauder estimates for nonuniformly elliptic problems by Mingione and the author in \cite{dm23a, dm25a}.
\subsection{Lavrentiev phenomenon}\label{lavlav} A classic of the Calculus of Variations, the Lavrentiev gap phenomenon occurs whenever variational integrals admit different infima if considered on the full class of admissible functions or on a smaller, but yet dense subclass of more regular functions, cf. \eqref{5.1}. This apparent anomaly was first observed by Lavrentiev \cite{lav26}. Later on, Manià \cite{man34}, considered the functional
$$
\mathcal{I}_{1}(w) :=\int_{0}^{1}(x-w^{3})^{2}\snr{w'}^{6}\dx,
$$
when defined
on the two different function spaces
\begin{flalign*}
&\mathcal{W}^{\infty}:=\left\{w\in W^{1,\infty}(0,1)\colon w(0)=0, \ w(1)=1\right\},\nonumber \\
&\mathcal{W}^{1}:=\left\{w\in W^{1,1}(0,1)\colon w(0)=0, \ w(1)=1\right\}\,.
\end{flalign*}
Obviously $u(x)=x^{1/3}$ minimizes $\mathcal{I}_{1}$ on $\mathcal{W}^{1}$, while explicit computations yield 
\begin{equation*}
\inf_{w\in \mathcal{W}^{\infty}}\mathcal{I}_{1}(w)>0=\inf_{w\in \mathcal{W}^{1}}\mathcal{I}_{1}(w)=\mathcal{I}_{1}(u)
\end{equation*}
cf. \cite[Lemma 4.43]{dac08}. Manià's example cannot be detected via standard finite element methods (i.e., taking piecewise affine functions that are Lipschitz continuous). One could argue that $\mathcal{I}_{1}$ lacks in coercivity. This issue was eventually fixed by Ball \& Mizel \cite{bm85} by looking at functional
$$
w \mapsto\int_{0}^{1}\varepsilon\snr{w'}^{2}+(x^{4}-w^{6})\snr{w'}^{m}\dx,
$$
with $m\ge 27$, and $\varepsilon>0$, that nonetheless presents the same pathological behavior as $\mathcal{I}_{1}$ in \eqref{5.1}. See also \cite{hm86} for the gap problem in stochastic control theory and  \cite{bk87} for numerical approximation schemes that spot lower energy singular minimizers despite the fact that the cost of any sequence in the admissible class of Lipschitz continuous functions is bounded away from the true minimum value\footnote{Note that in this case, when a minimizer over the smoother admissible class exists, the usual approximation schemes converge to this suboptimal solution.}. 

\subsection{Approximation in energy}\label{conyen}
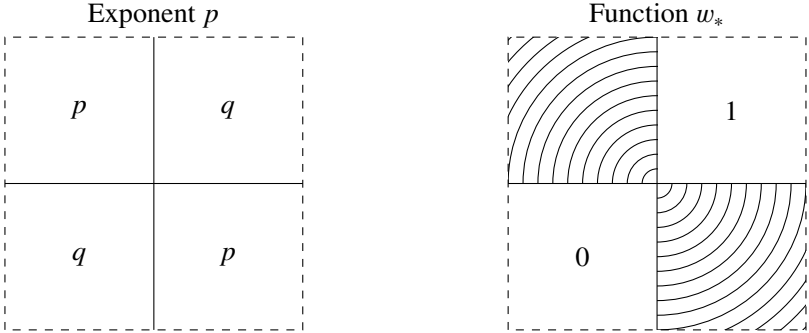
\begin{figure}[!ht]
  \centering
  \begin{tikzpicture}[scale=1.95]
    \node at (0,1.15) {Exponent~$p$};
    \draw[dashed] (-1,-1) -- (-1,+1) -- (+1,+1) -- (+1,-1) --cycle;
    \draw (-1,0) -- (1,0);
    \draw (0,-1) -- (0,1);
    \node at (0.5,0.5) {$q$};
    \node at (-0.5,-0.5) {$q$};
    \node at (-0.5,0.5) {$p$};
    \node at (0.5,-0.5) {$p$};
  \end{tikzpicture}
   \quad\quad\quad\quad\quad
  \qquad 
  \begin{tikzpicture}[scale=1.95]
    \node at (0,1.15) {Function $w_{*}$};
    \draw[dashed] (-1,-1) -- (-1,+1) -- (+1,+1) -- (+1,-1) --cycle;
    \draw (-1,0) -- (1,0);
    \draw (0,-1) -- (0,1);

    \node at (0.5,0.5) {$1$};
    \node at (-0.5,-0.5) {$0$};
    \begin{scope}
      \clip (-1,0) rectangle (0,1);
      \draw (0,0) circle (1.6);
      \draw (0,0) circle (1.5);
      \draw (0,0) circle (1.4);
      \draw (0,0) circle (1.3);
      \draw (0,0) circle (1.2);
      \draw (0,0) circle (1.1);
      \draw (0,0) circle (1.0);
      \draw (0,0) circle (0.9);
      \draw (0,0) circle (0.8);
      \draw (0,0) circle (0.7);
      \draw (0,0) circle (0.6);
      \draw (0,0) circle (0.5);
      \draw (0,0) circle (0.4);
      \draw (0,0) circle (0.3);
      \draw (0,0) circle (0.2);
      \draw (0,0) circle (0.1);
    \end{scope}
    \begin{scope}
      \clip (0,-1) rectangle (1,0);
      \draw (0,0) circle (1.6);
      \draw (0,0) circle (1.5);
      \draw (0,0) circle (1.4);
      \draw (0,0) circle (1.3);
      \draw (0,0) circle (1.2);
      \draw (0,0) circle (1.1);
      \draw (0,0) circle (1.0);
      \draw (0,0) circle (0.9);
      \draw (0,0) circle (0.8);
      \draw (0,0) circle (0.7);
      \draw (0,0) circle (0.6);
      \draw (0,0) circle (0.5);
      \draw (0,0) circle (0.4);
      \draw (0,0) circle (0.3);
      \draw (0,0) circle (0.2);
      \draw (0,0) circle (0.1);
    \end{scope}
  \end{tikzpicture}
   \caption{Zhikov’s checkerboard. Source \cite{bds20}.}
   \label{fig.1}
\end{figure}

The (absence of) Lavrentiev Phenomenon admits other, more application-friendly and general formulations beyond the one in \eqref{5.1}:
\begin{itemize}
\item Coincidence between the functional $\mathcal{F}$, and the LSM relaxation $\bar{\mathcal{F}}$.
\item Density of smooth maps in the related Lagrangian space\footnote{I.e., the class of all functions $w\in W^{1,1}(B)$ with finite energy $\mathcal{F}(w,\Omega)<\infty$.} \cite[Chapter 14]{jko94}.
\end{itemize}
Indeed, notice that whenever $w\in W^{1,p}(B)$ is such that $\bar{\mathcal{F}}(w,B)<\infty$ then 
there exists a sequence $
w_{i}\rightharpoonup w$   in $ W^{1,p}(B)$ such that $\mathcal{F}(w_{i},B)\to \bar{\mathcal{F}}(w,B)$. It follows that
$\mathcal{F}(w,B)=\bar{\mathcal{F}}(w,B)$ iff there exists a sequence \begin{equation}\label{conv}
w_{i}\rightharpoonup w \ \ \mbox{weakly in} \ \ W^{1,p}(B)\qquad \mbox{and}\qquad \mathcal{F}(w_{i},B)\to \mathcal{F}(w,B).
\end{equation}
This is called approximation in energy. If the integrand $F$ is nonautonomous, the occurrence of Lavrentiev phenomenon connects to the way space-depending coefficients  and gradient variable interact, in particular, \eqref{conv} may be violated even in simple cases. A paradigm of such irregularity phenomena is the celebrated Zhikov's checkerboard example \cite{zhi87, jko94}; consider the functional in \eqref{px} with $\Omega = \mathcal B_1\subset \mathbb{R}^{2}$  and, with $x\equiv (x_{1},x_{2})\in \mathbb{R}^{2}$, set
\begin{equation}\label{competi0}
\mf{p}(x):=\begin{cases}
\displaystyle
\, p\quad &\mbox{if} \ \ x_{1}x_{2}<0\vspace{.1mm}\\ \displaystyle
\, q\quad &\mbox{if} \ \ x_{1}x_{2}\geq 0,
\end{cases}\qquad \quad\quad  1<p<2<q\,.
\end{equation}
In polar coordinates $(x_{1},x_{2})=(\varrho \cos \theta,\varrho \sin \theta)$, define 
\begin{equation}\label{competi}
w_{*}(x):=\begin{cases}
    \displaystyle
    \, 1\quad &\mbox{if} \ \ x_{1},x_{2}>0\vspace{.1mm}\\
    \displaystyle
    \, \sin \theta\quad &\mbox{if} \ \ x_{2}>0>x_{1}\vspace{.1mm}\\
    \displaystyle
    \,  0\quad &\mbox{if} \ \ x_{1},x_{2}<0\vspace{.1mm}\\
    \displaystyle
    \,  \cos \theta \quad &\mbox{if} \ \ x_{1}>0>x_{2},
\end{cases}
\end{equation}
see Figure \ref{fig.1}. The Lagrangian space generated by the functional in \eqref{px} is given by
\begin{equation*}
W^{1,\mf{p}(\cdot)}(\mathcal B_1):=\left\{w\in W^{1,1}(\mathcal B_1)\colon \snr{Dw}^{\mf{p}(\cdot)}\in L^1(\mathcal B_1)\right\}.
\end{equation*}
The integral in  \eqref{px} is clearly well-defined and finite on $W^{1,\mf{p}(\cdot)}(\mathcal B_1)$. Zhikov proved that despite $w_{*}\in W^{1,\mf{p}(\cdot)}(\mathcal B_1)$, energy approximation by smooth maps in the sense of \eqref{conv} fails. This discussion highlights an interplay between gradient variable and space-depending coefficient, that is stronger than the one observed in the case of the $p$-Laplace energy: indeed now coefficients have no longer a perturbative role, and might dramatically affect the growth/ellipticity features of the integral\footnote{In fact, this already happens in the examples \eqref{dp}, \eqref{px} and \eqref{limiting}.}. We will deepen these facts in Section \ref{1.xx}, here we would simply note the huge amount of information obtained from the analysis of Lavrentiev phenomenon: functional analytic properties of Lagrangian spaces, relaxation of functionals, and, what matter the most here, possibility of constructing SOLAs\footnote{Solution Obtained by Limiting Approximations. These are solutions obtained as limits, in suitable sense, of solutions to more regular problems, like for instance minima of approximating uniformly elliptic functionals.}, i.e. building up energy approximations of general nonautonomous functionals, a crucial aspect of the analysis described in Section \ref{2.xx} (see comments after Theorem \ref{t2}). 

\subsection{Soft nonuniform ellipticity}\label{1.xx}

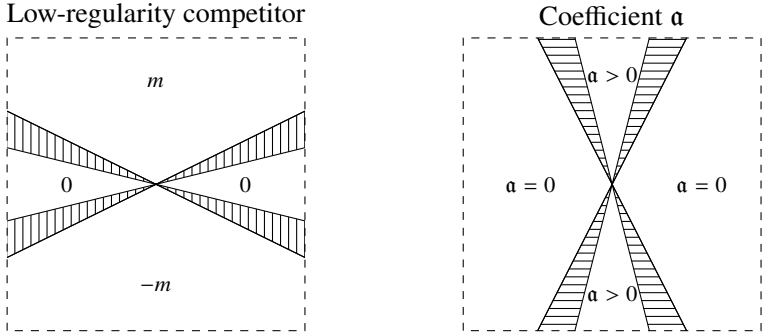
\begin{figure}[!ht]
  \centering
  \begin{tikzpicture}[scale=1.95]
    \node at (0,1.15) {Low-regularity competitor};
    \draw[dashed] (-1,-1) -- (-1,+1) -- (+1,+1) -- (+1,-1) --cycle;
    \node at (0,.7) {\scalebox{0.8}{$m$}};
    \node at (0,-.7) {\scalebox{0.8}{$-m$}};
    \draw (-1,-1/2) -- (1,1/2);
    \draw (-1,1/2)-- (1,-1/2);
    \node at (+.6,0) {\scalebox{0.8}{$0$}};
    \node at (-.6,0) {\scalebox{0.8}{$0$}};
    \filldraw[pattern=vertical lines] (-1,-1/2) -- (0,0) -- (-1,-1/4);
    \filldraw[pattern=vertical lines] (-1,+1/2) -- (0,0) -- (-1,+1/4);
    \filldraw[pattern=vertical lines] (+1,-1/2) -- (0,0) -- (+1,-1/4);
    \filldraw[pattern=vertical lines] (+1,+1/2) -- (0,0) -- (+1,+1/4);
  \end{tikzpicture}
  \quad\quad\quad\quad\quad
  \begin{tikzpicture}[scale=1.95]
    \node at (0,1.15) {Coefficient $\aaa$};
    \draw[dashed] (-1,-1) -- (-1,+1) -- (+1,+1) -- (+1,-1) --cycle;
    
    \node at (.6,0) {\scalebox{0.8}{$\aaa=0$}};
    \node at (-.55,0) {\scalebox{0.8}{$\aaa=0$}}; 
    \draw (-1/2,-1)--(+1/2,1);
    
    \draw(+1/2,-1)-- (-1/2,+1) ;
    \node at (0,.75) {\scalebox{0.8}{$\aaa>0$}};
    \node at (0,-.75) {\scalebox{0.8}{$\aaa>0$}};
    
    \filldraw[pattern=horizontal lines] (-1/2,-1) -- (0,0) -- (-1/4,-1);
    \filldraw[pattern=horizontal lines] (+1/2,-1) -- (0,0) -- (+1/4,-1);
    \filldraw[pattern=horizontal lines] (-1/2,+1) -- (0,0) -- (-1/4,+1);
    \filldraw[pattern =horizontal lines]
    (+1/2,+1) -- (0,0) -- (+1/4,+1);

  \end{tikzpicture}
  \caption{Low-regular competitor vs coefficient $\aaa$. Figure \ref{fig} is a modification of the one in \cite{bds20}.}
  \label{fig}
\end{figure}

Before focusing on classical, strong nonuniform ellipticity in the sense of \eqref{strong}, let us consider the softly nonuniformly elliptic case of the  double phase functional $\mathcal{D}$ in \eqref{dp}, where from now on we assume that 
\eqn{mainaaa}
$$0\leq \aaa\in C^{0,\alpha}_{\loc}(\Omega), \quad  \alpha \in (0,1]. $$  Behind a rather simple structure, energy $\mathcal{D}$ hides several interesting (ir)regularity aspects. In fact, building on previous two-dimensional constructions of Zhikov \cite{zhi95,zhi97} concerning the Lavrentiev phenomenon, it was discovered by Esposito \& Leonetti \& Mingione \cite{elm04} that, whenever exponents $(p,q)$ and $\alpha$ verify
\begin{equation}\label{4.0}
p<n<n+\alpha<q \qquad \left(\Longrightarrow \frac qp > 1+\frac \alpha n\right)
\end{equation}
it is possible to construct an integrand as in \eqref{dp} and \eqref{mainaaa}, and a Lipschitz-regular boundary datum $u_{0}\in W^{1,\infty}(\mathcal B_{1})$, such that the solution to the corresponding Dirichlet problem exhibits a one-point singularity that forbids its membership to $W^{1,q}_{\loc}(\mathcal B_{1})$, i.e., 
\begin{equation}\label{4.1}
u \mapsto \min_{w\in u_{0}+W^{1,p}_{0}(\mathcal B_{1})}\mathcal{D}(w,\mathcal B_{1}), \qquad u \not\in W^{1,q}_{\loc}(\mathcal B_{1})\,.
\end{equation} 
An even more striking example was shown by Fonseca \& Mal\'y \& Mingione \cite{fmm04}, who, under the same conditions of \cite{elm04} (but obviously with a different, much more complex coefficient $\aaa$), exhibited a solution to \eqref{4.1} whose set of its essential discontinuities forms a Cantor-type fractal of almost maximal Hausdorff dimension $n-p$\footnote{More precisely, in \cite{fmm04} for every $\varepsilon>0$ the authors construct an example where the  Hausdorff dimension of the set of essential discontinuity points of the minimizer is larger than $n-p-\varepsilon$, still keeping \eqref{4.0} and with  $q-p< \alpha +\varepsilon$. Note that any minimizer of $\mathcal{D}$  belongs to $W^{1,p}_{\loc}(\Omega)$ and that the Hausdorff dimension of the set of essential discontinuities of any $W^{1,p}_{\loc}$-function cannot exceed $n-p$.}. This result was eventually sharpened by Balci \& Diening \& Surnachev \cite{bds20,bds23}, who constructed bounded minima of $\mathcal{D}$ whose set of essential discontinuities has Hausdorff dimension equal to $n-p$, and therefore do not belong to any better Sobolev space than $W^{1,p}(B_{1})$. All in all, there are minima of convex, scalar, uniformly elliptic (in the classical sense) variational integrals that are as bad as any other competitor! This oddity confirms how delicate regularity issues are when dealing with nonuniform ellipticity, when if only a soft form of it is present. The key idea behind these examples is that the coefficient $\aaa$ vanishes where certain low-integrable competitor maps $w_*$ (no better than $W^{1,p}$) are non-constant, i.e., $\aaa(x)\snr{Dw_*(x)}^q=0$. These malicious competitors in \eqref{4.1} decrease the total energy below the energy of any possible $W^{1,q}$-competitor once highly varying traces ($m>0$ and $-m$) on the upper and lower parts of $\partial \mathcal B_{1}$ are chosen; see Figure \ref{fig}. 

This finally creates a non vanishing Lavrentiev gap and the situation is similar to the one described in Section \ref{1.xx}, see \eqref{competi0}-\eqref{competi}. Choosing a coefficient $\aaa$ with a fractal type geometry allows to distribute the corresponding singularities on a Cantor fractal with maximal Hausdorff dimension. In particular, the results in \cite{fmm04,bds23} highlight that, in presence of very mild types of nonuniformity as the one in \eqref{dp.1}$_1$, smooth coefficients do no longer guarantee regular solutions, cf. \eqref{3.0} and Remark \ref{rem1}. All such examples rely on the fact that under the considered condition the Lavrentiev gap of $\mathcal D$ does not vanish $ \mathcal{L}_{\mathcal D}\not= 0$. Condition \eqref{4.0} violates \eqref{4.2}, 
under which maximal gradient regularity is consequently expected. The first general regularity result is again in \cite{elm04}, where it is proved that if $u$ minimizes $\mathcal D$, then 
\eqn{capasso}
$$
\begin{cases}
\displaystyle \frac{q}{p} \leq 1 +\frac{\alpha}{n} & \Longrightarrow  \mathcal{L}_{\mathcal D}(\cdot,{\rm B})=0 \quad  \forall \ {\rm B}\Subset \Omega \\
& \Longrightarrow u\in W^{1,q+\delta}_{\loc}(\Omega;\mathbb R^n), \ \ \  0<\delta \equiv \delta(n,p,q, \alpha) \footnote{The additional equality case $q/p= 1+\alpha/n$ with respect to is a typical effect of the pointwise uniform ellipticity of the double phase integrand, allowing to get some extra small initial gradient integrability.}.
\end{cases}
$$
Then in \cite{cm15,bcm18} it is proved that 
\eqn{maximal}
$$
\frac{q}{p} \leq 1 +\frac{\alpha}{n} \Longrightarrow Du\in C^{1,\beta}_{\loc}(\Omega;\mathbb R^n), \ \ \ 0<\beta\equiv \beta(n,p,q, \alpha)<1.
$$
This is the best amount of regularity achievable in degenerate problems \cite{ura68}, already when $\aaa\equiv 0$. Let us highlight that the strategy in \cite{cm15,bcm18} is still perturbative. The soft nonuniformity of $\mathcal{D}$ in \eqref{dp.1} still allows to use certain perturbation methods of the type sketched in Section \ref{ues}, or, at least, "half" of them. More precisely, with reference to the "Freeze \& Compare" scheme in Section \ref{ues}, the main difficulty here relies in deriving homogeneous comparison estimates, cf. \eqref{3.6}, the validity of homogeneous reference estimates for suitable liftings being guaranteed by the pointwise uniform ellipticity in \eqref{dp.1}$_1$. Indeed, the proof of \eqref{maximal} conceptually goes as follows.
\begin{itemize} 
    \item Using pointwise uniform ellipticity \eqref{dp.1}$_1$. The "frozen" integrand with constant coefficient $\aaa_0$
    \begin{equation*}
    P_{0}(z):=\snr{z}^{p}+\aaa_0\snr{z}^{q},\qquad \quad \aaa_0\in [0, \infty)
    \end{equation*}
    is uniformly elliptic, and its ellipticity ratio is bounded independently of $\aaa_0$ by \eqref{dp.1}$_1$. With $B_{r}(x_{0})\Subset \Omega$ being an arbitrary  ball, this implies that the solution $v\in u+W^{1,p}_{0}(B_{r}(x_{0}))$ of the Dirichlet problem
    \eqn{solutio}
    $$
      v\mapsto \min_{u+W^{1,p}_{0}(B_{r}(x_{0}))} \int_{B_{r}(x_{0})}P_{0}(Dw)\dx,
$$
    features intrinsically homogeneous excess decay estimates of type \eqref{3.5}, i.e.:
    \begin{equation}\label{p0p0}
    \begin{cases}
   \displaystyle  \nr{P_{0}(Dv)}_{L^{\infty}(B_{r/2}(x_0))} \lesssim_{n,p,q} \mint_{B_{r}} P_{0}\left(Dv\right)\dx\\
    \displaystyle  \mint_{B_{\sigma}(x_0)}P_{0}\left(Dv-(Dv)_{B_{\sigma}(x_0)}\right)\dx\\
    \qquad \qquad  \displaystyle  \lesssim_{n,p,q} \left(\frac{\sigma}{\varrho}\right)^{\beta_1}\mint_{B_{\varrho}(x_0)} P_{0}\left( Dv-(Dv)_{B_{\varrho}(x_0)}\right)\dx 
     \end{cases}
    \end{equation} 
    for all concentric balls $B_{\sigma}(x_{0})\subset B_{\varrho}(x_{0})\Subset B_{r}(x_{0})$ and with $\beta_1\equiv \beta_1(n,p,q)\in (0,1)$ being independent of $\aaa_0$; compare with \eqref{3.5} to which \eqref{p0p0} reduce when $\aaa_0=0$. See \cite{lie91,dsv09} for a proof. Estimates \eqref{p0p0} hold for any choice of $1<p\leq q$; no upper bound on $q/p$ is needed. 
    \item Rebalancing the soft nonuniform ellipticity \eqref{dp.1}$_2$. This is the point where the optimal bound $q/p \leq 1+\alpha/n$ crucially enters the game. With $M$ large enough, to be quantitatively chosen, the $p$-phase is said to occur when\footnote{As usual $$[\aaa]_{0,\alpha;B_{r}(x_0)}=\sup_{x,y\in B_{r}(x_0); x\not=y} \, \frac{\snr{\aaa(x)-\aaa(y)}}{\snr{x-y}^\alpha}.$$}
    \eqn{pfase}
    $$
    \inf_{B_{r}(x_0)} \aaa \leq M [\aaa]_{0,\alpha;B_{r}(x_0)}r^\alpha\,.
    $$
    The $q$-phase instead occurs otherwise. In the $p$-phase the original minimizer is compared to the solution to \eqref{solutio} with $\aaa_0=0$. In the $(p,q)$-phase instead one chooses $\aaa_0=  \inf_{B_{r}(x_0)} \aaa $. Via delicate a comparison scheme based either on reverse H\"older inequalities \cite{cm15}, or on quantitative harmonic type approximation lemmas \cite{bcm18}, one arrives at 
    \eqn{compy}
    $$
        \mint_{B_{r/2}(x_{0})}P_{0}\left(Du-Dv\right)\dx\lesssim r^{\beta_{2}}\mint_{B_{r}(x_{0})}P_{0}\left(Du\right)\dx,
$$
    for some $\beta_{2}\in (0,1)$.  
    \item Estimates \eqref{p0p0} and \eqref{compy} are homogeneous estimates when considered with respect the intrinsic quantity $P_0(Du)$ and, as such, can now be matched and iterated in the setting of a delicate exit time argument based on the occurrence of \eqref{pfase} at various scales. This eventually yields gradient H\"older continuity for $u$. 
\end{itemize}
\noindent The viewpoint in \cite{cm15,bcm18} has been further developed, see \cite{bo20,bos22a,bgs22,bos22b,kl22,fsv24,bb25} and references therein for a non exhaustive list of relevant contributions. In particular, H\"asto \& Ok \cite{ho22a,ho22b,ho23} extended the outcome of \cite{bcm18} to a very large  family of softly nonuniformly elliptic problems, again, by means of a perturbative approach that strongly relies on the pointwise uniform ellipticity of the class of integrands considered. The whole body of techniques described above heavily relies on pointwise uniform ellipticity (soft nonuniform ellipticity) and fails to deliver results in presence of strong nonuniform ellipticity \eqref{strong}, even for the most basic model examples. 
\begin{remark} {\em It is worth noting that in the counterexample provided in \cite{fmm04}, the parameter $\alpha$ can be taken as any positive number $\alpha > 0$, and not necessarily restricted to the interval $\alpha \in (0,1]$, provided that $p$ and $q$ are chosen accordingly, as in \eqref{4.0}\footnote{In this case one takes $\aaa\in C^{[\alpha], \{\alpha\}}$.}. This shows that pathological minimizers can arise even when the coefficients have arbitrary degrees of smoothness, as long as $p$ and $q$ are sufficiently far apart.}
\end{remark}
\subsection{Nonuniformly elliptic Schauder estimates}\label{2.xx} 
The discussion in Section \ref{ues}  highlights that, while for PDEs driven by uniformly elliptic operators as in \eqref{unel} the classical Schauder implication
\begin{equation*}
\mbox{H\"older continuous coefficients} \ \Longrightarrow \ \mbox{H\"older continuous first derivatives}
\end{equation*} 
holds for energy solutions, already for the mild nonuniformity in \eqref{soft.1}-\eqref{soft.2} regularity of solutions is no longer guaranteed in general, not even by smooth ingredients. However, as seen in Section \ref{1.xx}, perturbation techniques are still available for soft nonuniformly elliptic equations or functionals, the main difficulty being getting close, in a homogeneous fashion, to a suitably "frozen" Dirichlet problem, whose solution features homogeneous excess decay estimates. Such a circle of ideas immediately breaks down already for the very simple model 
\begin{equation}\label{sm}
    \mathcal{F}_{a}(w,\Omega):=\int_{\Omega}a(x)F_{0}(Dw)\dx,\ \   1\lesssim a\in C^{0,\alpha}(\Omega), \ \alpha\in (0,1]
\end{equation}
where $F_{0}$ is as in \eqref{assf}. In fact, recalling \eqref{1.4}, in view of \eqref{assf} the best bound on the ellipticity ratio we have is still 
\begin{equation*}
\mathcal{R}(x,z)\approx\frac{\mbox{highest eigenvalue of} \ a(x)\partial_{zz}F_{0}(z)}{\mbox{lowest eigenvalue of} \ a(x)\partial_{zz}F_{0}(z)}\lesssim\snr{z}^{q-p} + 1\,.
\end{equation*}
Therefore any hope of recovering uniformly elliptic proof schemes as in Section \eqref{1.xx} vanishes. The validity of Schauder theory in the nonuniformly elliptic setting is a classical problem raised at several stages in the literature. For instance, in his MR review\footnote{Math. Rev. MR0749677} of Giaquinta \& Giusti's paper \cite{gg84}, Lieberman remarks how the authors' approach does not apply when uniform ellipticity conditions are violated, unless a priori Lipschitz bounds are assumed. However, these gradient estimates are only available either if equations or functionals are uniformly elliptic or if an unnatural amount of smoothness is imposed on coefficients and heavy regularity is a priori assumed on solutions. In other terms, regularity for nonuniformly elliptic problems with H\"older continuous coefficients could be delivered only for a priori Lipschitz continuous solutions\footnote{On the other hand in such cases the problems considered are essentially again uniformly elliptic as in the present setting the ellipticity ratio blows-up only when the gradient blows-up.}. The central role of gradient bounds is also remarked by Ivanov \cite[pages 7 and 15]{iva84}: in view of the results of Ladyzhenskaya \& Ural'tseva \cite{lu68}, the issue of solvability of boundary value problems for a nonuniformly elliptic or parabolic equation reduces to the possibility of constructing a priori estimates of the $L^{\infty}$-norm of the gradients of solutions for a suitable one-parameter family of similar equations. This means finding uniform a priori gradient estimates for regularized problems. The importance of deriving Lipschitz bounds in nonuniformly elliptic problems as a starting point for higher regularity was also stressed out by Giaquinta \& Giusti \cite[page 56]{gg84}. In addition to these considerations, Schauder theory does not always hold already in the case of soft nonuniform ellipticity \eqref{soft.1}-\eqref{soft.2}, cf. Section \ref{1.xx}. This emphasizes how delicate the issue of establishing Schauder theory is for strongly nonuniformly elliptic problems \eqref{strong}. 
Such a longstanding problem has been settled in \cite{dm23a,dm25a}, through a novel approach based on a direct proof of gradient boundedness, that does not make any use of H\"older gradient estimates\footnote{The same approach also leads to settle the longstading issue of proving gradient H\"older continuity of minima of nondifferentiable, nonuniformly elliptic variational integrals. The additional difficulty in this case is that in general the Euler-Lagrange equation is not available. See \cite{dm23a}.}. This alters the classical paradigm described in the third point of Remark \ref{rem1}. The following results appear in \cite{dm25a}. 
\begin{theorem}\label{t1}
    Under assumptions \eqref{assf} and \eqref{4.2}, let $u\in W^{1,p}_{\loc}(\Omega)$ be a local minimizer of the LSM-extension $\bar{\mathcal{F}}$ of the functional $\mathcal{F}$ in \eqref{fun}. Then 
    \begin{itemize}
    \item $Du$ is locally H\"older continuous in $\Omega$.
    \item $u$ is also a minimizer of $\mathcal{F}$. 
    \end{itemize}
    \end{theorem}
    Recalling \eqref{newminima} we have 
\begin{cor}\label{t1.2}
    Under assumptions \eqref{assf} and \eqref{4.2}, let $u\in W^{1,p}_{\loc}(\Omega)$ be a local minimizer of the functional $\mathcal{F}$ in \eqref{fun}. If $\mathcal{L}_{\mathcal F}(u,{\rm B})=0$ holds for all balls ${\rm B}\Subset \Omega$, then $u$ is a minimizer of $\bar{\mathcal{F}}$ and hence $Du$ is locally H\"older continuous in $\Omega$. \end{cor}
\noindent If integrand $F$ is nondegenerate or nonsingular, in the sense that $\mu>0$ in \eqref{assf}, we recover the classical Schauder type implication \eqref{3.0}.
\begin{theorem}\label{t2}
In Theorem \ref{t1} and Corollary \ref{t1.2}, assume that $\mu>0$ and $\alpha\in (0,1)$ in \eqref{assf}, and that $\partial_{zz}F$ is continuous. Then $Du\in C^{0,\alpha}_{\loc}(\Omega;\mathbb{R}^{n})$.
\end{theorem}
\noindent In particular, Theorem \ref{t1} and Corollary \ref{t2} highlight that the LSM-extension $\bar{\mathcal{F}}$ provides a selection principle capable of identifying/producing regular minimizers of the original functional $\mathcal{F}$, while excluding those that are irregular due to the appearance of Lavrentiev type phenomena. In Corollary \ref{t1.2} the role of the assumption $\mathcal{L}_{\mathcal F}(u,{\rm B})=0$ is rather natural: via the occurrence of the approximation in energy it is possible to approximate the original minimizer with a sequence of minimizers of more regular functionals for which a priori regularity estimates, in fact described in the next section, apply. The Lavrentiev gap functional vanishes in a large variety of situations, for instance detailed in \cite[Section 5]{elm04}. As an example, a condition of the type $F(x,z)\approx G(z)$, for $\snr{z}$ large, where $G\colon \mathbb R^n \to [0,\infty)$ is a convex function, guarantees that $\mathcal L_{\mathcal F}\equiv 0$. As a consequence we have
\begin{cor} \label{main}
Let $u\in W^{1,1} (\Omega)$ be a minimizer of the functional in \eqref{sm}
where 
$F_0$ satisfies \eqref{assf}$_{1-4}$ and \eqref{4.2} holds. Then, $Du$ is locally H\"older continuous in $\Omega$. If $\mu>0$, then $u\in C^{1,\alpha}_{\loc}(\Omega)$ when $\alpha <1$.
 \end{cor}
 As a matter of fact, an interesting twist is that on several occasions assumption \eqref{4.2} itself guarantees that the Lavrentiev gap vanishes. For instance, this happens when 
$F(x,z)\approx \snr{z}^p+\aaa(x)\snr{z}^q$ 
holds for large $\snr{z}$ and \eqref{4.2} is satisfied. This last fact does not even require that $F$ satisfies \eqref{assf} and holds for general Carath\'eodory integrands; compare with \eqref{capasso}.  
\subsection{Key steps in the proof of Theorem \ref{t1}}\label{ultima} The proof offers a general body of methods and perspectives that might be useful in a number of problems with lack of ellipticity and/or regularity of external ingredients. The main arguments leading to Theorem \ref{t1} can be outlined as follows\footnote{For simplicity, we describe how a priori estimates are obtained, results then come by approximation procedures for which we refer to \cite{dm25a}. Keep in mind the comments made after Theorem \ref{t2}.}.
\begin{itemize}
    \item \textbf{Almost Lipschitz continuity via fractional Moser's iteration.} The gradient of minima belongs to $L^{t}_{\loc}$ for all $1\le t<\infty$ via a combined use of interpolation inequalities and a fractional Moser type iteration. The first step is showing that a relaxed bound on $q/p$, i.e., 
    \begin{equation}\label{5.3}
 q < p+ 
\begin{cases}
\,  p\alpha  /2& \mbox{if $p<2$}\\
\,  \alpha & \mbox{if $2 \leq p \leq n$}\\
\,  p\alpha/n & \mbox{if $n < p$}\,.
\end{cases}
    \end{equation}
implies that any (locally) bounded minimizer has locally integrable gradient with arbitrary large, finite power (local boundedness of minima under the assumptions used here follows by \cite{hs21}). In this respect, \eqref{5.3} is optimal by counterexamples working for bounded minimizers \cite{elm04}, where minima fail to be $W^{1,q}$-regular; compare \eqref{5.3} with \eqref{4.0}. In the autonomous case\footnote{with the exception of the specific double phase case, see \cite{bcm18}.}, results of this type were known, and classical, starting by the work of Ladyzhenskaya \& Ural'tseva \cite{lu70}, see also \cite{uu84,cho92,ckp11,bcm18, ddp24}, where dimension-free gap bounds are proved to be effective when in presence of bounded minimizers. The techniques used for them  require differentiating the Euler-Lagrange equation \eqref{1.7}, which is impossible in presence of H\"older coefficients. In \cite{dm25a} higher integrability estimates are instead achieved by means of fine Besov spaces techniques, sometimes employed in uniformly elliptic problems with a certain lack of ellipticity, such as degenerate PDEs in the Heisenberg group \cite{dom04}, or in nonlocal problems \cite{bl17,bls18,gl23,dkln24,gl24,bdlms24}. The underlying principle is that dealing with strong nonuniform ellipticity requires similar efforts as those needed to compensate the lack of strong forms of monotonicity. This is delicate to handle and heavily relies on a sharp numerology to keep the optimal conditions \eqref{5.3}. Specifically, we prove that
    \begin{equation}\label{5.4}
        \left\|\frac{\tau_{h}^{2}u}{\snr{h}^{1+\varepsilon_{i}}}\right\|_{L^{t_{i}}} \le c_i
    \end{equation}
holds for $h \in \mathbb{R}^n$, $|h|>0$ sufficiently small, and with integro-differentiability exponents $t_{i} \nearrow \infty$, $0< \varepsilon_i  \searrow 0$ and $c_i\to \infty$. In \eqref{5.4}, $\tau^{2}_{h}$ denotes the double finite differences operator in the direction $h$\footnote{With $\texttt{t}>0$, $h \in \mathbb{R}^n$, set $\Omega_{\texttt{t}}:=\left\{x\in \Omega\colon \dist(x,\partial \Omega)>\texttt{t}\right\}$; the finite difference operators $\tau_{h}\colon L^{1}(\Omega;\mathbb{R}^{n})\mapsto L^{1}(\Omega_{|h|};\mathbb{R}^{n})$, $\tau^{2}_{h}\colon L^{1}(\Omega;\mathbb{R}^{n})\mapsto L^{1}(\Omega_{2|h|};\mathbb{R}^{n})$ are defined as $ \tau_{h}w(x):=w(x+h)-w(x)$ and  $\tau_{h}^{2}w:=\tau_{h}(\tau_{h}w)$, respectively.}. Via basic embedding properties of Besov spaces, \eqref{5.4} implies that $Du\in L^{t_{i}}_{\loc}$ for every $i$, and therefore in all Lebesgue spaces with finite exponent, locally. We stress that the validity of \eqref{5.3} guarantees that sequence $\{t_{i}\}_{i\in \mathbb{N}}$ in \eqref{5.4} diverges. Estimates leading to \eqref{5.4} implement a fractional Moser type iteration that does not directly lead to $Du \in L^{\infty}_{\loc}$, as the integrability gain from $t_{i}$ to $t_{i+1}$ is linear rather than geometric. A fractional Moser's iteration is a natural idea to use in this setting as it relies on the basic observation that $\alpha$-H\"older continuity means almost $\alpha$-fractional differentiability, in the sense that 
$$
x\mapsto \partial_{z} F(x,z)\stackrel{\eqref{assf}_{4}}{\in} C^{0,\alpha} \ \Longrightarrow \ x\mapsto \partial_{z} F(x,z)\in W^{s,t} ,
$$
for all $s\in (0,\alpha)$, $t\in [1,\infty)$, thus aligning with classical Moser's iteration in the differentiable setting, where in fact one uses that $x\mapsto \partial_{z} F(x,z)\in W^{1,\infty}$.  
\item \textbf{Fractional Caccioppoli inequalities on level sets.} A standard formulation of gradient Caccioppoli inequalities on level sets for autonomous, quasilinear uniformly elliptic PDEs like for instance $\textnormal{div}(\snr{Du}^{p-2} Du)=0$, reads as
\begin{equation}\label{5.5}
    \varrho^{2}\mint_{B_{\varrho/2}(x_{0})}\snr{D(\snr{Du}^{p}-\kappa)_{+}}^{2}\lesssim \mint_{B_{\varrho}(x_{0})}(\snr{Du}^{p}-\kappa)_{+}^{2}\dx,
\end{equation}
for any $\kappa\ge 0$ and all balls $B_{\varrho}(x_{0})\Subset \Omega$. Estimate \eqref{5.5} builds on classical Bernstein method, that prescribes differentiating the equation to prove that suitable convex functions of $\snr{Du}$, such as for instance $\snr{Du}^{p}$, are subsolutions to linear, uniformly elliptic equations with measurable coefficients, hence, they satisfy \eqref{5.5}. In this respect see the seminal work of Uhlenbeck \cite{uhl77}. This yields local boundedness for $\snr{Du}$ since \eqref{5.5} grants its membership to upper De Giorgi classes, i.e., classical De Giorgi's iteration applies. This argument breaks down in presence of merely H\"older continuous coefficients, as the Euler-Lagrange equation is no longer differentiable. Nonetheless the fractional, renormalized Caccioppoli inequality on level sets
\begin{flalign}\label{5.7}
     \varrho^{2\beta-n} [(|Du|^p-\kappa)_{+}]_{\beta,2;B_{\varrho/2}(x_{0})}^2\lesssim& M^{2\textnormal{\texttt{b}}} \mint_{B_{\varrho}(x_{0})}(|Du|^p-\kappa)_{+}^{2}\dx\nonumber \\
     &+ M^{2\textnormal{\texttt{b}}}\varrho^{2\alpha}\mint_{B_{\varrho}(x_{0})}1+\snr{Du}^{m}\dx, 
\end{flalign}
 holds whenever $M\ge \nr{Du}_{L^{\infty}(B_{\varrho}(x_{0}))}$ for some $\beta\in (0,1)$, $m\in (1,\infty)$, on any ball $B_{\varrho}(x_{0})\Subset \Omega$. In \eqref{5.7} the number $\texttt{b}\equiv \texttt{b}(n,p,q,\alpha)$ is an explicit function such that $\texttt{b}\to 0$ when $q\to p$.  There are three major differences between the classical \eqref{5.5} and the formulation \eqref{5.7}, as introduced in \cite{dm23a, dm25a}\footnote{Fractional Caccioppoli type inequalities in the setting of nonlinear potential theory were pioneered by Mingione \cite{min07,min11}, where also fractional De Giorgi's classes are used. See also \cite{ccv11,dkp16,coz17} for other instances of fractional De Giorgi classes in the nonlocal setting.}. First, the fractional derivatives of  $(\snr{Du}^{p}-\kappa)_{+}$   appearing in \eqref{5.7} replace the full derivatives on the left-hand side of \eqref{5.5}, which are no longer available due to the nondifferentiability of the coefficients. Second, the third term  in \eqref{5.7} is a "size"  term that captures the contribution of the coefficients and requires careful handling throughout the analysis. Third,  
 the renormalization constant $M$  is introduced to make the first two terms in  \eqref{5.7} appear homogeneous in $\snr{Du}^{p}$, mirroring the structure of \eqref{5.5}, and thereby enabling iteration. This strategy also draws on an approach developed by Beck \& Mingione \cite{bm20}, which prescribes—when working with inequalities like \eqref{5.5} in the nonuniformly elliptic setting—to exploit the fact that constants deteriorate only polynomially during De Giorgi iterations used to prove $L^\infty$-estimates\footnote{This behavior is less widely known, as constants deteriorate exponentially instead when proving Hölder estimates, which are the ultimate aim of the De Giorgi–Nash–Moser theory.}. Accordingly, careful tracking of the quantity $M^{2\texttt{b}}$ throughout the iteration process raises the possibility of ultimately reabsorbing it. In this respect, exponent $\texttt{b}$ accounts for the power type growth of the ellipticity ratio, while $\beta\equiv \beta(\alpha)\in (0,1)$ records the rate of H\"older continuity of coefficients. The sharp threshold \eqref{4.2} plays a pivotal role in linking these analytic features, ultimately yielding local Lipschitz regularity for minimizers. Crucially, achieving this under the optimal threshold \eqref{4.2} requires that \eqref{5.7} is employed in its sharpest possible form, ensuring that no loss of information occurs at any stage of the derivation. The construction of the renormalized, fractional De Giorgi classes in \eqref{5.7} is based on a Littlewood-Paley type argument traceable to the dyadic decomposition of Besov functions \cite{tri01}, but implemented at a nonlinear level. It was first introduced by Kristensen \& Mingione \cite{km05,km06} in the study of singular sets of minima of multiple integrals. In the present setting the usual atoms appearing in dyadic decompositions are replaced by solutions to nonlinear, autonomous problems, whose regularity is quantified via suitable improved forms of related a priori estimates due to Bella \& Schäffner \cite{bs20,bs24,sch24}. This grants the sharp value of exponent $\texttt{b}$ in \eqref{5.7}, that allows to preserve the optimal threshold in \eqref{4.2}.
\item \textbf{Nonlinear potentials.} A key aspect of the De Giorgi type iteration described before lays in its potential theoretic nature. In fact, the size information encoded in the third  term of \eqref{5.7} is carried along iterations by a general class of nonlinear Havin-Maz'ya-Wolff type potentials \cite{hm72} of the form
\begin{equation*}
{\bf P}_{\sigma}^{\vartheta}(f;x_0,r) := \int_0^r \varrho^{\sigma} \left(  \mint_{B_{\varrho}(x_0)} \snr{f} \dx \right)^{\vartheta} \frac{\textnormal{d}\varrho}{\varrho},
\end{equation*}
with $\sigma>0$, $\vartheta\ge 0$ and $f\in L^{1}(B_{r}(x_{0}))$. For specific choices of $\sigma$ and $\vartheta$, potential $\mathbf{P}^{\vartheta}_{\sigma}$ becomes classical Riesz and Wolff potentials, that are usually employed when investigating fine properties of solutions to linear or nonlinear, nonhomogeneous PDEs. More precisely, potentials act as "ghosts" of the representation formula\footnote{Representation formulae via fundamental solutions are obviously not available in nonlinear problems.}, granting sharp pointwise bounds for solutions and their gradients under minimal conditions on external data, see \cite{km94,km12,km13,km18}. This means that potentials have been mainly used to deal with certain borderline regularity conditions on external ingredients. In \cite{dm23a,dm23b,dm25a,ddp24} instead, potentials are employed as carriers of size information along iterations, and eventually apply their mapping properties \cite[Section 4]{dm23a} to control the resulting quantities. Also at this stage, in absence of gradient higher integrability, condition \eqref{4.2} is crucial to ensure that potentials embed in the right Lebesgue spaces, \cite{dm23a,dm23b}.
\end{itemize}
The approach sketched above yields local Lipschitz continuity for minima. In turn, gradient boundedness makes the nonuniformity of functional \eqref{fun} (and of its LSM-relaxation \eqref{LSM}) immaterial, so gradient H\"older continuity is achieved by a delicate adaptation of more standard perturbation methods \cite[Section 10]{dm23a}, \cite[Section 5]{dm25a}. See also the third point in Remark \ref{rem1}.


\begin{ack}
The author thanks Professor Anna Balci for sharing the drawings in Figures \ref{fig.1}-\ref{fig}.
\end{ack}

\begin{funding}
This work was supported by the University of Parma through the action "Bando di Ateneo 2024 per la ricerca".
\end{funding}


\end{document}